\documentclass[usletter,11pt]{article}

\usepackage{fullpage}

\usepackage[latin1]{inputenc}
\usepackage{graphicx}

\usepackage{amsmath}
\usepackage{amsthm}

\usepackage{amssymb}
\usepackage{bbm}
\usepackage{stmaryrd}
\usepackage{comment}
\usepackage{ifthen}

\usepackage{amsthm}

\newtheorem{lemma}{Lemma}
\newtheorem{corollary}{Corollary}
\newtheorem{theorem}{Theorem}
\newtheorem{proposition}{Proposition}
\newtheorem{remark}{Remark}

\theoremstyle{definition}
\newtheorem{definition}{Definition}

\newcommand{\mtiret}{-\hspace{-2mm}-}

\newcommand{\infM}{<_{\mathfrak{m}}}

\newcommand{\M}{{\mathfrak{m}}}

\newcommand{\barM}{{\overline{\mathfrak{m}}}}
\newcommand{\infbarM}{<_{\overline{\mathfrak{m}}}}

\newcommand{\infMii}{<_{\mathfrak{m}_{i+1}}}

\newcommand{\Cat}{\mathrm{Cat}}

\title{A new combinatorial identity for unicellular maps, via a direct bijective approach}

\author{
   Guillaume Chapuy\\
\small  Department of Mathematics,
   Simon Fraser University,
   Burnaby, B.C. V5A 1S6
\thanks{Supported by a CNRS/PIMS postdoctoral fellowship. This work was mainly done while the author was a PhD student in LIX, Laboratoire d'Informatique de l'\'Ecole Polytechnique (France). The author acknowledges support from the grant ERC StG 208471 - ExploreMap.}
}

\begin{document}

\maketitle

\begin{abstract}
A unicellular map, or one-face map, is a graph embedded in an orientable surface such that its complement is a topological disk. 
In this paper, we give a new viewpoint to the structure of these objects, by describing a decomposition of any unicellular map into a unicellular map of smaller genus.
 This gives a new combinatorial identity for the number $\epsilon_g(n)$ of unicellular maps of size~$n$ and genus~$g$. Contrarily to the Harer-Zagier recurrence formula, this identity is recursive in only one parameter (the genus).

Iterating the construction gives an explicit bijection between unicellular maps and plane trees with distinguished vertices, which gives a combinatorial explanation (and proof) of the fact that $\epsilon_g(n)$ is the product of the $n$-th Catalan number by a polynomial in $n$. The combinatorial interpretation also gives a new and simple formula for this polynomial. Variants of the problem are considered, like bipartite unicellular maps, or unicellular maps with cubic vertices only.


{\small
\begin{tabbing}
{\bf Keywords:}~~\= Polygon gluings, Bijection, Harer-Zagier numbers.
\end{tabbing}
}
\end{abstract}


\section{Introduction.}

A \emph{unicellular map} is a graph embedded in a compact orientable surface, in such a way that its complement is a topological polygon.
Equivalently, a unicellular map can be viewed as a polygon, with an even number of edges, in which edges have been pasted pairwise in order to create a closed orientable surface. The number of handles of this surface is called the \emph{genus} of the map.

These objects are reminiscent in combinatorics, and have been considered in many different contexts. The numbers of unicellular maps of given size and genus appear in random matrix theory as the moments of the Gaussian Unitary Ensemble~(see~\cite{LaZv}). In the study of characters of the symmetric group, unicellular maps appear as factorisations of cyclic permutations~\cite{Jackson:countingcycles}. According to the context, unicellular maps are also called \emph{one-face maps}, \emph{polygon gluings}, or \emph{one-border ribbon graphs}. Sometimes, their duals, \emph{one-vertex maps}, are considered.
The most famous example of unicellular maps are \emph{plane} unicellular maps, which, from Jordan's lemma, are exactly plane trees, enumerated by the Catalan numbers.

The first result in the enumeration of unicellular maps in positive genus was obtained by Lehman and Walsh~\cite{Lehman-Walsh-genus-I}. Using a direct recursive method, relying on formal power series, they expressed the number $\epsilon_g(n)$ of unicellular maps with $n$ edges on a surface of genus $g$ as follows:
\begin{eqnarray}
\label{eq:LW}
\epsilon_g(n)= \sum_{\gamma \vdash g } \frac{(n+1)\dots(n+2-2g-l(\gamma))}{2^{2g} \prod_i c_i! (2i+1)^{c_i}} \textrm{Cat}(n), 
\end{eqnarray}
where the sum is taken over partitions $\gamma$ of $g$, $c_i$ is the number of parts $i$ in $\gamma$, $l(\gamma)$ is the total number of parts, and $\textrm{Cat}(n)$ is the $n$-th Catalan number. This formula has been extended by other authors (\cite{Goupil-Schaeffer}).

Later, Harer and Zagier~\cite{Harer-Zagier}, via matrix integrals techniques, obtained the two following equations, known respectively as the Harer-Zagier \emph{recurrence} and the Harer-Zagier \emph{formula}:
\begin{eqnarray}\label{eq:HZrecurrence}
(n+1) \epsilon_g(n) = 2(2n-1) \epsilon_g(n-1) + (2n-1)(n-1)(2n-3) \epsilon_{g-1}(n-2), \\
 \sum_{g\geq 0} \epsilon_{g}(n) y ^{n+1-2g} = \frac{(2n)!}{2^nn!}\sum_{i\geq 1} 2^{i-1} {n \choose i-1} {y \choose i}.
\label{eq:HZformula}
\end{eqnarray}
Formula~\ref{eq:HZformula} has been retrieved by several authors, by various techniques. A combinatorial interpretation of this formula was given by Lass~\cite{Lass}, and the first bijective proof was given by Goulden and Nica~\cite{Goulden-Nica}. Generalizations were given for bicolored, or multicolored maps \cite{Jackson:countingcycles,Schaeffer-Vassilieva}.

The purpose of this paper is to give a new angle of attack to the enumeration of unicellular maps, at a level which is much more combinatorial than what existed before. Indeed, until now no bijective proof (or combinatorial interpretation) of Formulas~\ref{eq:LW} and~\ref{eq:HZrecurrence} are known. As for Formula~\ref{eq:HZformula}, it is concerned with some \emph{generating polynomial} of the numbers $\epsilon_g(n)$: 
in combinatorial terms, the bijections in \cite{Goulden-Nica,Schaeffer-Vassilieva} concern maps which are \emph{weighted} according to their genus, by an additional coloring of their vertices, but the genus does not appear explicitely in the constructions. For example, one cannot use these bijections to sample maps of given genus and size.

On the contrary, this article is concerned with the structure of unicellular maps themselves, at fixed genus. We investigate in details the way the unique face of such a map interwines with itself in order to create the handles of the surface. We show that, in each unicellular map of genus $g$, there are $2g$ special "places", which we call \emph{trisections}, that concentrate, in some sense, the handles of the surface. Each of these places can be used to \emph{slice} the map to a unicellular map of lower genus.
Conversely, we show that a unicellular map of genus $g$ can always be obtained in $2g$ different ways by gluing vertices together in a map of lower genus. 
In terms of formulas, this leads us to the new combinatorial identity:
\begin{eqnarray}\label{eq:identity}
 2g  \cdot \epsilon_g(n)\!\!\!\!\! &=& \!\!\!\!\! {n+3-2g\choose 3} \epsilon_{g-1}(n) + {n+5-2g\choose 5} \epsilon_{g-2}(n) + \dots + {n+1\choose 2g+1} \epsilon_{0}(n) \\
&=& \!\!\!\!\! \sum_{p=0}^{g-1} {n+1-2p \choose 2g-2p+1} \epsilon_p(n).
\end{eqnarray}
The main advantage of this identity is that it is recursive \emph{only in the genus}: the size $n$ is fixed. For a given $g$, this enables one to compute directly the formula giving $\epsilon_g(n)$, by iteration. From the combinatorial viewpoint, this enables one to construct maps of fixed genus and size very easily.

When iterated, our bijection shows that all unicellular maps can be obtained in a canonical way from plane trees by successive gluings of vertices, 
hence giving the first explanation to the fact that $\epsilon_g(n)$ is the product of a polynomial $R_g(n)$ by the $n$-th Catalan number. More precisely, we obtain the formula $\epsilon_g(n)=R_g(n)\mathrm{Cat}(n)$ with:
\begin{eqnarray}\label{eq:Rg}
R_g(n) = \sum_{0=g_0<g_1<\dots<g_r=g}\ \prod_{i=1}^r \frac{1}{2g_i}{{n+1-2g_{i-1}}\choose{2(g_i-g_{i-1})+1}},
\end{eqnarray}
which comes with a clear combinatorial interpretation. This interpretation, and the one of certain properties of the polynomial $R_g$, answers questions of Zagier~\cite[p159]{LaZv}.\\

\noindent {\bf Asymptotic case.} In the paper~\cite{Chapuy:trisections}, we presented a less powerful bijection, that worked only for an asymptotically dominating subset of all unicellular maps. The bijection presented here is really a generalization of the bijection of~\cite{Chapuy:trisections}, in the sense that it coincides with it when specialized to those dominating maps. However, new difficulties and structures appear in the general case, and there is an important gap between the combinatorial results in~\cite{Chapuy:trisections} and the ones of this paper. 

\noindent {\bf Extended abstract.} A extended abstract of this paper was presented at the conference FPSAC'09 (Austria, July 2009). 

\noindent {\bf Acknowledgements.} I am indebted to Olivier Bernardi, and to Gilles Schaeffer, for very stimulating discussions. Thanks also to Emmanuel Guitter for allowing me to mention his computation in Section~\ref{subsec:guitter}.

\section{Unicellular maps.} 

\label{sec:maps}

\subsection{Permutations and ribbon graphs.}
\begin{figure}[h]
\centerline{\includegraphics[scale=1.1]{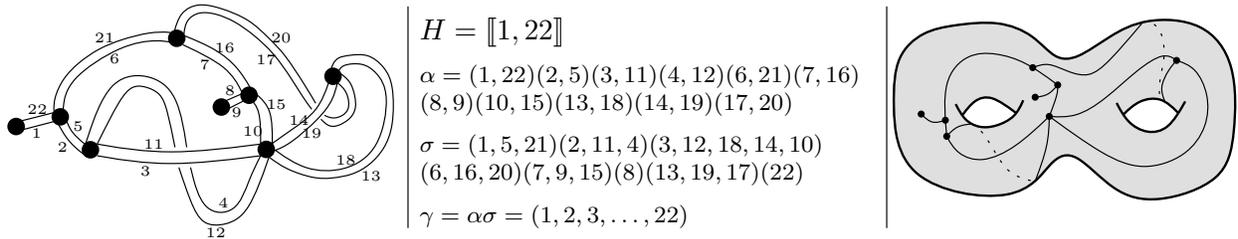}}
\caption{A unicellular map with $11$ edges, $8$ vertices, and genus $2$: (a) ribbon graph; (b) permutations; (c) topological embedding.}
\label{fig:map}
\end{figure}

Rather than talking about topological embeddings of graphs, we work with a combinatorial definition of unicellular maps:
\begin{definition}
A \emph{unicellular map} $\M$ of size $n$ is a triple $\M=(H,\alpha,\sigma)$, where $H$ is a set of cardinality $2n$, $\alpha$ is an involution of $H$ without fixed points, 
and $\sigma$ is a permutation of $H$ such that $\gamma=\alpha\sigma$ has only one cycle.
The elements of $H$ are called the \emph{half-edges} of $\M$. The cycles of $\alpha$ and $\sigma$ are called the \emph{edges} and the \emph{vertices} of $\M$, respectively, and the permutation $\gamma$ is called the \emph{face} of $\M$.
\end{definition}
Given a unicellular map $\M=(H,\sigma,\alpha)$, its associated (multi)graph $G$ is the graph whose edges are given by the cycles of $\alpha$, vertices by the cycles of $\sigma$, and the natural incidence relation $v\sim e$ if $v$ and $e$ share an element. Moreover, we draw each edge of $G$ as a \emph{ribbon}, where each side of the ribbon represents one half-edge; we decide which half-edge corresponds to which side of the ribbon by the convention that, if a half-edge $h$ belongs to a cycle $e$ of $\alpha$ and $v$ of $\sigma$, then $h$ is the right-hand side of the ribbon corresponding to $e$, when considered entering $v$. Furthermore, we draw the graph $G$ in such a way that around each vertex $v$, the \emph{counterclockwise ordering} of the half-edges belonging to the cycle $v$ is given by that cycle: we obtain a graphical object called the \emph{ribbon graph} associated to $\M$, as in Figure~\ref{fig:map}(a). Observe that the unique cycle of the permutation $\gamma=\alpha\sigma$ interprets as the sequence of half-edges visited when \emph{making the tour} of the graph, keeping the graph on its left.
%

A \emph{rooted} unicellular map is a unicellular map carrying a distinguished half-edge $r$, called the root. These maps are considered up to relabellings of $H$ preserving the root, i.e. two rooted unicellular maps $\M$ and $\M'$ are considered the same if there exists a permutation $\pi: H \rightarrow H'$, such that $\pi(r)=r'$, $\alpha=\pi^{-1}\alpha'\pi$, and $\sigma=\pi^{-1}\sigma'\pi$. In this paper, \emph{all} unicellular maps will be rooted, even if not stated.

Given a unicellular map $\M$ of root $r$ and face $\gamma=\alpha\sigma$, we define the linear order $<_{\scriptscriptstyle \M}$ on $H$ by setting:
$$r <_{\scriptscriptstyle \M} \gamma(r) <_{\scriptscriptstyle \M} \gamma_2(r) <_{\scriptscriptstyle \M} \dots <_{\scriptscriptstyle \M} \gamma^{2n-1}(r).
$$
In other words, if we relabel the half-edge set $H$ by elements of $\llbracket 1, 2n \rrbracket$ in such a way that the root is $1$ and the tour of the face is given by the permutation $(1,\dots,2n)$, the order $<_{\scriptscriptstyle \M}$ is the natural order on the integers. However, since in this article we are going to consider maps with a fixed half-edge set, but a changing permutation $\gamma$, it is more convenient (and prudent) to define the order $<_{\scriptscriptstyle \M}$ in this way.

Unicellular maps can also be interpreted as graphs embedded in a topological surface, in such a way that the complement of the graph is a topological polygon. If considered up to homeomorphism, and suitably rooted, these objects are in bijection with ribbon graphs. See~\cite{Mohar-Thomassen}, or the example of Figure~\ref{fig:map}(c).
The \emph{genus} of a unicellular map is the genus, or number of handles, of the corresponding surface. If a unicellular map of genus $g$ has $n$ edges and $v$ vertices, then Euler's characteristic formula says that $v=n+1-2g.$
From a combinatorial point of view, this last equation can also be taken as a definition of the genus.

\subsection{The gluing operation.}


\begin{figure}[h]
\centerline{\includegraphics[scale=0.8]{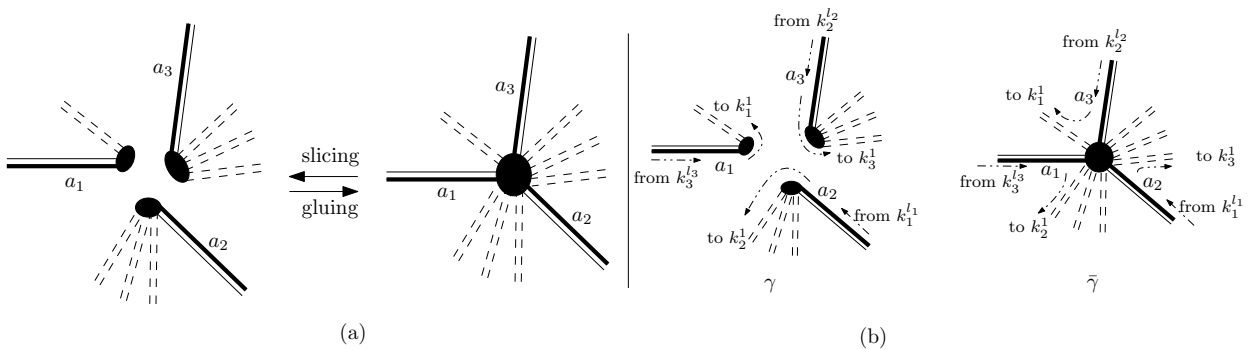}}
\caption{(a) The gluing and slicing operations. (b) The "proof" of Lemma~\ref{lemma:facegluing}.}
\label{fig:gluing}
\end{figure}

We let $\M=(H,\alpha,\sigma)$ be a unicellular map of genus $g$, and 
 $a_1<_{\scriptscriptstyle \M} a_2 <_{\scriptscriptstyle \M} a_3$ be three half-edges of $\M$ \emph{belonging to three distinct vertices}. Each half-edge $a_i$ belongs to some vertex $v_i=(a_i,h^1_i,\dots h^{m_i}_i)$, for some $m_i\geq 0$. We define the permutation 
$$
 \bar v:= (a_1,h^1_2,\dots h^{m_2}_2,a_2,h^1_3,\dots h^{m_3}_3,a_3,h^1_1,\dots h^{m_1}_1),
$$ 
and we let $\bar\sigma$ be the permutation of $H$ obtained by deleting the cycles $v_1$, $v_2$, and $v_3$, and replacing them by $\bar v$. 
The transformation mapping $\sigma$ to $\bar \sigma$ interprets combinatorially as the \emph{gluing} of the three half-edges $a_1,a_2,a_3$, as shown on Figure~\ref{fig:gluing}(a). We have:
\begin{lemma}\label{lemma:facegluing}
The map 
$\barM:=(H,\alpha,\bar\sigma)$ is a unicellular map of genus $g+1$. 
Moreover, if we let 
$$\gamma=\alpha\sigma=(a_1, k^1_1,\dots k^{l_1}_1,a_2, k^1_2,\dots k^{l_2}_2,a_3, k^1_3,\dots k^{l_3}_3)$$
 be the face permutation of $\M$, then the face premutation of $\barM$ is  given by:
$$
\bar \gamma = (a_1, k^1_2,\dots k^{l_2}_2,a_3, k^1_1,\dots k^{l_1}_1,a_2, k^1_3,\dots k^{l_3}_3)
$$
\end{lemma}
\begin{proof}
In order to prove that $\M$ is a well-defined unicellular map, it suffices to check that its face is given by the long cycle $\bar\gamma$ given in the lemma. This is very easy by observing that the only half-edges whose image is not the same by $\gamma$ and by $\bar \gamma$ are the three half-edges $a_1,a_2,a_3$, and that by construction $\bar \gamma(a_i)=\alpha\bar\sigma(a_i)=\alpha\sigma(a_{i+1})=\gamma(a_{i+1})$. For a more "visual" explanation, see Figure~\ref{fig:gluing}(b).

Now, by construction, $\barM$ has two less vertices than $\M$, and the same number of edges, so from Euler's formula it has genus $g+1$ (intuitively, the gluing operation has created a new "handle").
\end{proof}


\subsection{Some intertwining hidden there, and the slicing operation.}

\begin{figure}
\centerline{\includegraphics[scale=1.2]{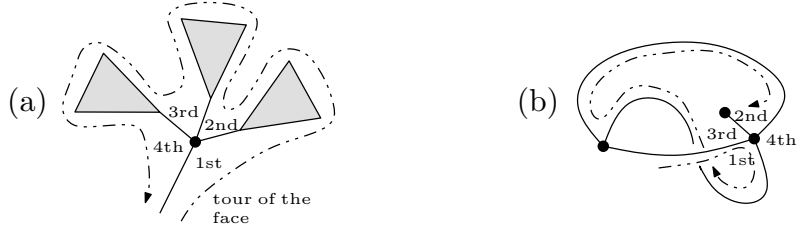}}
\caption{(a) In a plane tree, the tour of the face always visits the half-edges around one vertex in \emph{counterclockwise} order; (b) in positive genus (here in genus $1$), things can be different.}
\label{fig:intertwining}
\end{figure}

The aim of this paper is to show that all unicellular maps of genus $g+1$ can be obtained \emph{in some canonical way} from unicellular maps of genus $g$ from the operation above. This needs to be able to "revert" (in some sense) the gluing operation, hence to be able to determine, given a map of genus $g+1$, which vertices may be "good candidates" to be sliced-back to a map of lower genus.

Observe that in the unicellular map $\barM$ obtained after the gluing operation, the three half-edges $a_1$, $a_2$, $a_3$ appear \emph{in that order} around the vertex $\bar v$, whereas they appear in the \emph{inverse} order in the face $\bar\gamma$. Observe also that this is very different from what we observe in the planar case: if one makes the tour of a plane tree, with the tree on its left, then one necessarily visits the different half-edges around each vertex in \emph{counterclockwise order} (see Figure~\ref{fig:intertwining}). 
Informally, one could hope that, in a map of positive genus, those places where the vertex-order does not coincide with the face-order hide some "intertwining" (some handle) of the map, and that they may be used to slice-back the map to lower genus.

We now describe the \emph{slicing} operation, which is nothing but the gluing operation, taken at reverse.
We let $\barM=(H,\alpha,\bar\sigma)$ be a map of genus $g+1$, and three half-edges $a_1, a_2, a_3$ belonging to a same vertex $\bar v$ of~$\barM$. We say that $a_1, a_2, a_3$ are \emph{intertwined} if they do not appear in the same order in $\bar \gamma=\alpha\bar\sigma$ and in $\bar \sigma$. In this case, we write 
$\bar v=(a_1,h^1_2,\dots h^{m_2}_2,a_2,h^1_3,\dots h^{m_3}_3,a_3,h^1_1,\dots h^{m_1}_1)$, and we let $\sigma$ be the permutation of $H$ obtained from $\bar \sigma$ by replacing the cycle $\bar v$ by the product $(a_1, h_1^1, \dots h_{m_1}^1)(a_2, h_1^2, \dots h_{m_2}^2)(a_3, h_1^3, \dots h_{m_3}^3).$
\begin{lemma}\label{lemma:faceslicing}
The map $\M=(H,\alpha,\sigma)$ is a well-defined unicellular map of genus $g$. If we let 
$$\bar \gamma = (a_1, k^1_2,\dots k^{l_2}_2,a_3, k^1_1,\dots k^{l_1}_1,a_2, k^1_3,\dots k^{l_3}_3)$$
 be the unique face of $\barM$, then the unique face of $\M$ is given by: 
$$\gamma=\alpha\sigma=(a_1, k^1_1,\dots k^{l_1}_1,a_2, k^1_2,\dots k^{l_2}_2,a_3, k^1_3,\dots k^{l_3}_3).$$

The gluing and slicing operations are inverse one to the other.
\end{lemma}
\begin{proof}
The proof is the same as in Lemma~\ref{lemma:facegluing}: it is sufficient to check the expression given for $\gamma$ in terms of~$\bar \gamma$, which is easily done by checking the images of $a_1,a_2,a_3$.
\end{proof}

\subsection{Around one vertex: up-steps, down-steps, and trisections.}

\begin{figure}
 \begin{minipage}{.56\linewidth}
\vspace{4.5mm}
\centerline{\includegraphics[scale=1.0]{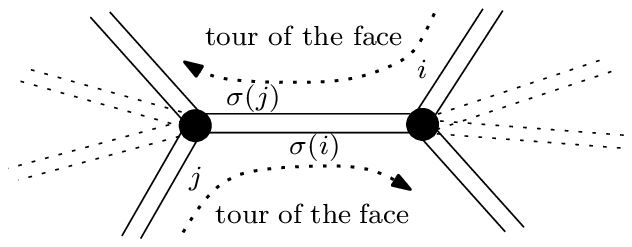}}
\vspace{3mm}
\caption{The main argument in the proof of the trisection lemma: the tour of the face visits $i$ before $\sigma(i)$ if and only if it visits $\sigma(j)$ before $j$, unless $
\sigma(i)$ or $\sigma(j)$ is the root of the map.}
\label{fig:trisectionproof}
 \end{minipage}\hfill 
\begin{minipage}{.38\linewidth}
\centerline{\includegraphics[scale=1.0]{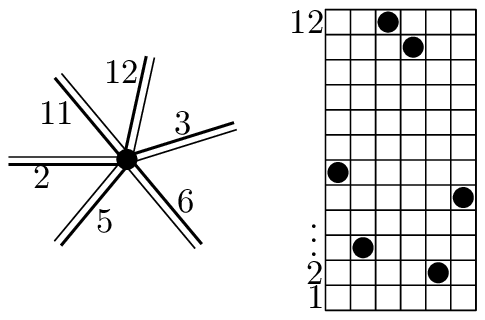}}
\caption{A vertex  $(6,3,12,11,2,5)$ in a map with $12$-half-edges, and its diagram representation (the marked half-edge is~$6$).}
\label{fig:diagram}
 \end{minipage} 
\end{figure}

Let $\M=(H,\alpha,\sigma)$ be a map of face permutation $\gamma=\alpha\sigma$. For each vertex $v$ of $\M$, we let $\min_{\scriptscriptstyle \M}(v)$ be the minimal half-edge belonging to $v$, for the order $<_{\scriptscriptstyle \M}$. Equivalently, $\min_{\scriptscriptstyle \M} (v)$ is the first half-edge from which one reaches $v$ when making the tour of the map, starting from the root. Given a half-edge $h\in H$, we note $V(h)$ the unique vertex it belongs to (i.e. the cycle of $\sigma$ containing it).
\begin{definition}
We say that a half-edge $h\in H$ is an \emph{up-step} if $h<_{\scriptscriptstyle \M} \sigma(h)$, and that it is a \emph{down-step} if $\sigma(h)\leq_{\scriptscriptstyle \M} h$. A down-step $h$ is called a \emph{trisection} if $\sigma(h)\neq \min_{\scriptscriptstyle \M} V(h)$, i.e. if $\sigma(h)$ is not the minimum half-edge inside its vertex.
\end{definition}
As illustrated on Figure~\ref{fig:intertwining}, trisections are specific to the non-planar case (there are no trisections in a plane tree), and one could hope that trisections "hide", in some sense, the handles of the surface. Before making this more precise, we state the following lemma, which is the cornerstone of this paper:
\begin{lemma}[The trisection lemma]\label{lemma:trisection}
Let $\M$ be a unicellular map of genus $g$. Then $\M$ has exactly $2g$ trisections.
\end{lemma}
\begin{proof}
We let $\M=(H,\alpha,\sigma)$, and $\gamma=\alpha\sigma$. We let $n_+$ and $n_-$ denote the number of up-steps and down-steps in $\M$, respectively. Then, we have $n_-+n_+=2n$, where $n$ is the number of edges of $\M$. Now, let $i$ be a half-edge of $\M$, and $j=\sigma^{-1}\alpha\sigma(i)$. Observe that we have $\sigma(j)=\gamma(i)$, and $\gamma(j)=\sigma(i)$. Graphically, $i$ and $j$ lie in two "opposite" corners of the same edge, as shown on Figure~\ref{fig:trisectionproof}. On the picture, it seems clear that if the tour of the map visits $i$ before $\sigma(i)$, then it necessarily visits $\sigma(j)$ before $j$ (except if the root is one of these four half-edges) so that, roughly, there must be almost the same number of up-steps and down-steps. More precisely, let us distinguish three cases.

First, assume that $i$ is an up-step. Then we have $i<_{\scriptscriptstyle \M}\sigma(i)=\gamma(j)$. Now, by definition of the total order~$<_{\scriptscriptstyle \M}$, $i <_{\scriptscriptstyle \M} \gamma(j)$ implies that $\gamma(i)\leq_{\scriptscriptstyle \M}\gamma(j)$.
Hence, $\sigma(j)\leq_{\scriptscriptstyle \M} \gamma(j)$, which, by definition of $<_{\scriptscriptstyle \M}$ again, implies that $\sigma(j)\leq_{\scriptscriptstyle \M} j$ (here, we have used that $\sigma(j)\neq \gamma(j) $ since $\alpha$ has no fixed point). Hence, if $i$ is an up-step, then $j$ is a down-step.

Second, assume that $i$ is a down-step, \emph{and} that $\gamma(j)$ is not equal to the root of $\M$. In this case, we have $j<_{\scriptscriptstyle \M}\gamma(j)$, and $\gamma(j)=\sigma(i)\leq_{\scriptscriptstyle \M} i=\sigma(j)$. Hence $j<_{\scriptscriptstyle \M} \sigma(j)$, and $j$ is an up-step.

The third and last case is when $i$ is a down step, \emph{and} $\gamma(j)$ is the root $r$ of $\M$. In this case, $j$ is the maximum element of $H$ for the order $<_{\scriptscriptstyle \M}$, so that it is necessarily a down-step. 

Therefore we have proved that each edge of $\M$ (more precisely, each cycle of $\sigma^{-1}\alpha\sigma$) is associated to one up-step and one down-step, except a special one that has two down-steps.
Consequently, there are exactly two more down-steps that up-steps in the map $\M$, i.e.: $n_-=n_++2$. Recalling that $n_-+n_+=2n$, this gives $n_-=n+1$.

Finally, each vertex of $\M$ carries exactly one down-step which is not a trisection (its minimal half-edge). Hence, the total number of trisections equals $n_--v$, where $v$ is the number of vertices of $\M$. Since from Euler's characteristic formula, $v$ equals $n+1-2g$, the lemma is proved.
\end{proof}

\section{Making the gluing operation injective.}

\label{sec:diagrams}

We have defined above an operation that glues a triple of half-edges, and increases the genus of a map. In this section, we explain that, if we restrict to \emph{certain types} of triples of half-edges, this operation can be made reversible.

\subsection{A diagram representation of vertices.}
\label{subsec:diagrams}

We first describe a graphical visualisation of vertices. Let $v$ be a vertex of $\M$, with a distinguished half-edge $h$. We write $v=(u_0,u_1,\dots,u_m)$, with $u_0=h$. We now consider a grid with $m+1$ columns and $2n$ rows. Each row represents an element of $H$, and the rows are ordered from the bottom to the top by the total order $<_{\scriptscriptstyle \M}$ (for example the lowest row represents the root). Now, for each $i$, inside the $i$-th column, we plot a point at the height corresponding to the half-edge $u_i$. We say that the obtained diagram is \emph{the diagram representation of $v$, starting from $h$}. 
In other words, if we identify $\llbracket 1,2n\rrbracket$ with $H$ via the order $<_{\scriptscriptstyle \M}$, the diagram representation of $v$ is the graphical representation of the sequence of labels appearing around the vertex $v$.
If one changes the distinguished half-edge $h$, the diagram representation of $v$ is changed by a circular permutation of its columns.
Figure~\ref{fig:diagram} gives an example of such a diagram (where the permutation $\gamma$ is in the form $\gamma=(1,2,3,\dots)$).

\begin{figure}[h]
\centerline{\includegraphics[scale=0.85]{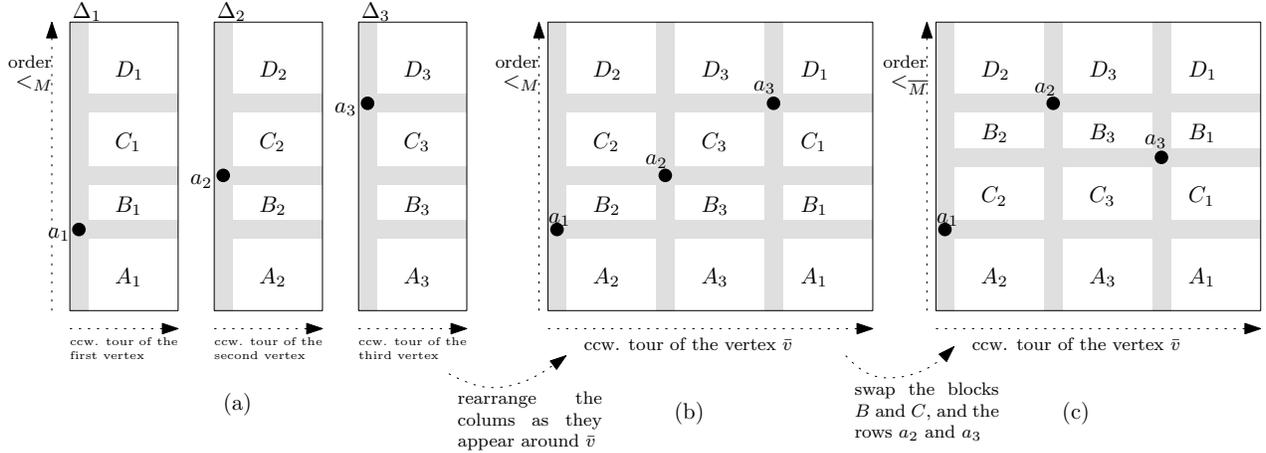}}
\caption{The gluing operation visualized on diagrams. (a) the diagrams before gluing; (b) a temporary diagram, where we the columns represent the counterclockwise turn around $\bar v$, but the rows still represent the original permutation $\gamma$; (c) the final diagram of the new vertex in the new map, where the rows represent the permutation $\bar\gamma$.}
\label{fig:diaggluing}
\end{figure}

The gluing operation is easily visualised on diagrams. We let as before $a_1<_{\scriptscriptstyle \M} a_2<_{\scriptscriptstyle \M} a_3$ be three half-edges belonging to distinct vertices in a unicellular map $\M$, and we let $\Delta_1,\Delta_2,\Delta_3$ be their corresponding diagrams. We now consider the three horizontal rows corresponding to $a_1$, $a_2$, and $a_3$: they separate each diagram $\Delta_i$ into four blocks (some of which may be empty). We give a name to each of these blocks: $A_i,B_i,C_i,D_i$, from bottom to top, as on Figure~\ref{fig:diaggluing}(a).

We now juxtapose $\Delta_2,\Delta_3,\Delta_1$ together, from left to right,  and we rearrange the three columns containing $a_1,a_2,a_3$ so that these half-edges appear in that order: we obtain a new diagram (Figure~\ref{fig:diaggluing}(b)), whose columns represent the order of the half-edges around the vertex $\bar v$. But the \emph{rows} of that diagram are still ordered according to the order $<_{\scriptscriptstyle \M}$. 
In order to obtain the diagram representing $\bar v$ in the new map $\barM$, we have to rearrange the rows according to $<_{\scriptscriptstyle \barM}$.
We let $A$ be the union of the three blocks $A_i$ (and similarly, we define $B$, $C$, and $D$). We know that the face permutation of $\M$ has the form
$
{\gamma=(\mtiret A \mtiret,a_1,\mtiret B \mtiret,a_2,\mtiret C \mtiret,a_3,\mtiret D \mtiret) },
$
where by $\mtiret A \mtiret$, we mean "all the elements of $A$, appearing in a certain order". Now, from the expression of $\bar \gamma$ given in Lemma~\ref{lemma:facegluing},  the permutation $\bar\gamma$ is:
$
\bar \gamma=(\mtiret A \mtiret,a_1,\mtiret C \mtiret,a_3,\mtiret B \mtiret,a_2,\mtiret D \mtiret),
$
where inside each block, the half-edges appear in the same order as in $\gamma$. In terms of diagrams, this means that the diagram representing $\bar v$ in the new map $\barM$ can be obtained by \emph{swapping} the block $B$ with the block $C$, and the row corresponding to $a_2$ with the one corresponding to $a_3$: see Figure~\ref{fig:diaggluing}(c). To sum up, we have:
\begin{lemma}
The diagram of the vertex $\bar v$ in the map $\barM$ is obtained from the three diagrams $\Delta_1,\Delta_2,\Delta_3$ by the following operations, as represented on Figure~\ref{fig:diaggluing}:\\
\noindent - Juxtapose $\Delta_2,\Delta_3,\Delta_1$ (in that order), and rearrange the columns containing $a_1, a_2,a_3$, so that they appear in that order from left to right. \\
\noindent - Exchange the blocks $B$ and  $C$, and swap the rows containing $a_2$ and $a_3$.
\end{lemma}
Observe that, when taken at reverse, Figure~\ref{fig:diaggluing} gives the way to obtain the diagrams of the three vertices resulting from the \emph{slicing} operation of three intertwined half-edges $a_1,a_2,a_3$ in the map $\barM$.

\begin{remark}\label{rem:ordre}
The slicing operation does not change the order $\infbarM$ for half-edges which appear strictly between the root and the minimum half-edge of the three vertices $\{a_1,a_2,a_3\}$. Precisely if $w_1\infbarM w_2 \infbarM \dots \infbarM w_r$ are elements of $H$ such that $w_r \infbarM a_2$, then Lemma~\ref{lemma:faceslicing} (or, more visually, Figure~\ref{fig:diaggluing}) implies that we have :
$$
w_1 \infM w_2 \infM \dots \infM w_r \infM a_1
$$
in the map $\M$. The reverse statement is also true.
\end{remark}

\subsection{Gluing three vertices: trisections of type {\bf I}.}

In this section, we let $v_1, v_2, v_3$ be three distinct vertices in the map $\M$. We let $a_i:=\min_{\scriptscriptstyle \M} v_i$, and, up to re-arranging the three vertices, we may assume (and we do) that $a_1<_{\scriptscriptstyle \M}a_2<_{\scriptscriptstyle \M}a_3$. We let $\Delta_1$, $\Delta_2$, $\Delta_3$ be the three corresponding diagrams. Since in each diagram the marked edge is the minimum in its vertex, observe that the blocks $A_1$, $A_2$, $B_2$, $A_3$, $B_3$, $C_3$ do not contain any point. We say that they are \emph{empty}, and we note: $A_1=A_2=B_2=A_3=B_3=C_3=\varnothing$.

We now glue the three half-edges $a_1,a_2,a_3$ in $\M$: we obtain a new unicellular map $\barM$, with a new vertex $\bar v$ resulting from the gluing.
Now, let $\tau$ be the element preceding $a_3$ around $\bar v$ in the map $\barM$. Since $A_3=B_3=C_3=\varnothing$, we have either $\tau\in D_3$ or $\tau=a_2$, so that in both case $a_3<_{\scriptscriptstyle \barM} \tau$. 
Moreover, $a_3$	 in not the minimum inside its vertex (the minimum is $a_1$). Hence, $\tau$ is a \emph{trisection} of the map $\barM$. We let $\Phi(\M,v_1,v_2,v_3)=(\barM,\tau)$ be the pair formed by the new map $\barM$ and the newly created trisection $\tau$.

\begin{figure}[h]
\centerline{\includegraphics[scale=0.8]{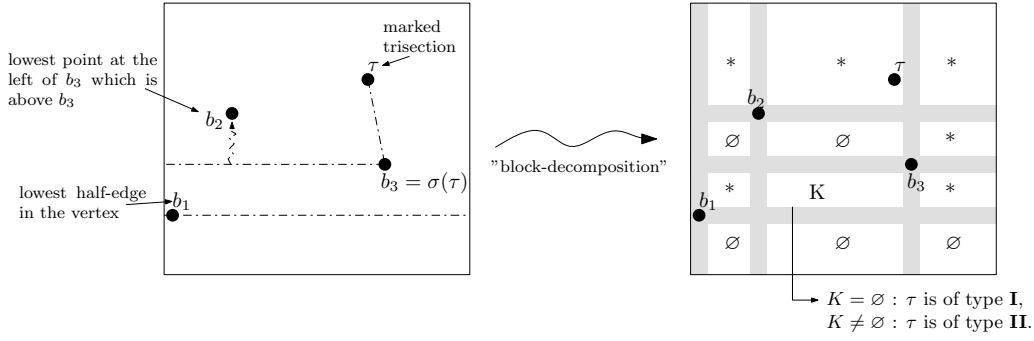}}
\caption{Trisections of type {\bf I} and {\bf II}.}
\label{fig:trisectionmarquee}
\end{figure}

It is clear that given $(\barM,\tau)$, we can inverse the gluing operation. Indeed, it is easy to recover the three half-edges $a_1$ (the minimum of the vertex), $a_3$ (the one that follows $\tau$), and $a_2$ (observe that, since $B_2$ and $B_3$ are empty, $a_2$ is the smallest half-edge on the left of $a_3$ which is greater than $a_3$). Once $a_1,a_2,a_3$ are recovered, it is easy to recover the map $\M$ by \emph{slicing} $\bar v$ at those three half-edges. This gives:
\begin{lemma}
The mapping $\Phi$, defined on the set of unicellular maps with three distinguished (unordered) vertices, is injective.
\end{lemma}
It is natural to ask for the image of $\Phi$: in particular, can we obtain all pairs $(\barM,\tau)$ in this way~? The answer needs the following definition (see Figure~\ref{fig:trisectionmarquee}):
\begin{definition}\label{def:trisection}
Let $\barM=(H,\alpha,\bar \sigma)$ be a map of genus $g+1$, and $\tau$ be a trisection of $\barM$. We let $\bar v=V(\tau)$, $b_1=\min_{\scriptscriptstyle \barM}(\bar v)$, and we let $\Delta$ be the diagram representation of $\bar v$, starting from the half-edge $b_1$. We let $b_3=\sigma(\tau)$ be the half-edge following $\tau$ around $\bar v$, and we let $b_2$ be the minimum half-edge among those which appear before $b_3$ around $\bar v$ and which are greater than $b_3$ for the order~$<_{\scriptscriptstyle \barM}$.
 
The rows and columns containing $b_1,b_2,b_3$ split the diagram $\Delta$ into twelve blocks, five of which are necessarily empty, as in Figure~\ref{fig:trisectionmarquee}.
We let $K$ be second-from-left and second-from-bottom block.
We say that $\tau$ is \emph{a trisection of type {\bf I}} is $K$ is empty, and that $\tau$ is \emph{a trisection of type {\bf II}} otherwise.
\end{definition}
The following proposition is the half way to our main result:
\begin{proposition}\label{prop:type1}
The mapping $\Phi$ is a bijection between the set $\mathcal{U}^3_{g}(n)$ of unicellular maps of genus $g$ with $n$ edges and three distinguished vertices,
and the set $\mathcal{D}^I_{g+1}(n)$ of unicellular maps of genus $g+1$ with $n$ edges and a distinguished trisection of type {\bf I}. 
\end{proposition}
\begin{proof}
We already know that $\Phi$ is injective.


We let $\M$ be a unicellular map of genus $g$ with three distinguished vertices $v_1, v_2, v_3$, and $\barM$ be the map obtained, as above, by the gluing of $M$ by the half-edges $a_1=\min_{\scriptscriptstyle \M} v_1$, $a_2=\min_{\scriptscriptstyle \M} v_2$, $a_3=\min_{\scriptscriptstyle \M} v_3$ 
(we assume again that $a_1<_{\scriptscriptstyle \M}a_2<_{\scriptscriptstyle \M}a_3$). We let $\bar \Delta$ be the diagram representation of the new vertex $\bar v$ obtained from the gluing in the map $\barM$, and we use the same notations for the blocks as in Section~\ref{subsec:diagrams}. We also let $\tau=\sigma^{-1}(a_3)$ be the created trisection, and we use the notations of Definition~\ref{def:trisection} with respect to the trisection $\tau$, so that $b_3=a_3$. Then, since $a_1=\min_{\scriptscriptstyle \barM} \bar v$, we have $a_1=b_1$, and since the blocks $B_2, B_3$, are empty, we have $b_2=a_2$. Hence, the block $C_3$ of Figure~\ref{fig:diaggluing}(c) coincides with the block $K$ of Figure~\ref{fig:trisectionmarquee}. Since $C_3$ is empty, $\tau$ is a trisection of type {\bf I}. Therefore the image of $\Phi$ is included in $\mathcal{D}^I_{g+1}(n)$.

Conversely, let $\barM=(H,\alpha,\bar\sigma)$ be a map of genus $g+1$, and $\tau$ be a trisection of type {\bf I} in $\barM$. We let $b_1, b_2, b_3$ and $K$ be as in Definition~\ref{def:trisection}. First, since $b_1<_{\scriptscriptstyle \barM} b_3<_{\scriptscriptstyle \barM} b_2$, these half-edges are intertwined, and we know that the slicing of $\barM$ by these half-edges creates a well-defined unicellular map $\M$ of genus $g$ (Lemma~\ref{lemma:faceslicing}). Now, if we compare Figures~\ref{fig:trisectionmarquee} and~\ref{fig:diaggluing}, we see that the result of the slicing is a triple of vertices $v_1,v_2,v_3$, such that each half-edge $b_i$ is the minimum in the vertex $v_i$: indeed, the blocks $A_1,A_2,A_3,B_2,B_3$ are empty by construction, and the block $C_3=K$ is empty since $\tau$ is a trisection of type {\bf I}. Hence we have $\Phi(\M,v_1,v_2,v_3)=(\barM, \tau)$, so that the image of $\Phi$ exactly equals the set $\mathcal{D}^I_{g+1}(n)$.
\end{proof}

\subsection{Trisections of type {\bf II}.}

Of course, it would be nice to have a similar result for trisections of type {\bf II}.
Let $\barM=(H,\alpha,\bar\sigma)$ be a map of genus $g+1$ with a distinguished trisection $\tau$ of type {\bf II}. We let $b_1,b_2,b_3$ and $K$ be as in Definition~\ref{def:trisection} and Figure~\ref{fig:trisectionmarquee}, and we let $\M$ be the result of the slicing of $\barM$ at the three half-edges $b_1,b_2,b_3$. If we use the notations of Figure~\ref{fig:diaggluing}, with $a_i=b_i$, we see that we obtain three vertices, of diagrams $\Delta_1,\Delta_2,\Delta_3$, such that $A_1=A_2=B_2=A_3=B_3=\varnothing$. Hence, we know that $a_1=\min_{\scriptscriptstyle \M}(v_1)$, that $a_2=\min_{\scriptscriptstyle \M}(v_2)$, and that $a_2<\min_{\scriptscriptstyle \M}(v_3)$. Observe that, contrarily to what happened in the previous section, the block $C_3=K$ is \emph{not} empty, therefore $a_3$ is \emph{not} the minimum inside its vertex.

Now, we claim that $\tau$ \emph{is still a trisection} in the map $\M$. Indeed, by construction, we know that $\tau$ belongs to $D_3$ (since, by definition of a trisection, it must be above $a_3$ in the map $\barM$, and since $B_3$ is empty). Hence we still have $a_3 <_{\scriptscriptstyle \M} \tau$ in the map $\M$. Moreover, we have clearly $\sigma(\tau)=a_3$ in $\M$ (since $\tau$ is the rightmost point in the blocks $C_3\cup D_3$), and it follows that $\tau$ is a trisection in $\M$.

We let $\Gamma(\barM,\tau)=(\M,v_1,v_2,\tau)$ be the $4$-tuple consisting of the new map $\M$, the two first vertices $v_1$ and $v_2$ obtained from the slicing, and the trisection $\tau$. It is clear that $\Gamma$ is injective: given $(\barM,v_1,v_2,\tau)$, one can reconstruct the map $\barM$ by letting $a_1=\min v_1$, $a_2=\min v_2$, and $a_3=\sigma(\tau)$, and by gluing back together the three half-edges $a_1,a_2,a_3$. Conversely, we define:
\begin{definition}
We let $\mathcal{V}_g(n)$ be the set of $4$-tuples $(\M,v_1,v_2,\tau)$, where 
$\M$ is a unicellular map of genus $g$ with $n$ edges, and where $v_1$, $v_2$, and $\tau$ are respectively two vertices and a trisection of $M$ such that: 
\begin{eqnarray}\label{eq:contrainteordre}
\min_{\scriptscriptstyle \M} v_1 <_{\scriptscriptstyle \M} \min_{\scriptscriptstyle \M} v_2 <_{\scriptscriptstyle \M} \min_{\scriptscriptstyle \M} V(\tau).
\end{eqnarray}
Given $(\M,v_1,v_2,\tau)\in\mathcal{V}_G(n)$, we let $\barM$ be the map obtained from the gluing of the three half-edges $\min v_1$, $\min v_2$, and $\sigma(\tau)$, and we let $\Psi(\M,v_1,v_2,\tau):=(\barM,\tau)$.
\end{definition}
We can now state the following proposition, that completes Proposition~\ref{prop:type1}:
\begin{proposition}\label{prop:type2}
The mapping $\Psi$ is a bijection between
the set $\mathcal{V}_g(n)$ of unicellular maps of genus $g$ with $n$ edges and a distinguished triple $(v_1,v_2,\tau)$ satisfying Equation~\ref{eq:contrainteordre}, 
and the set $\mathcal{D}^{II}_{g+1}(n)$ of unicellular maps of genus $g+1$ with $n$ edges and a distinguished trisection of type {\bf II}.
\end{proposition}
\begin{proof}
In the discussion above, we have already given a mapping $\Gamma: \mathcal{D}^{II}_{g+1}(n) \rightarrow \mathcal{V}_g(n)$, such that $\Psi\circ\Gamma$ is the identity on $\mathcal{D}_{g+1}^{II}(n)$.

Conversely, let $(\M,v_1,v_2,\tau)\in\mathcal{V}_{g}(n)$, and let $a_1=\min v_1$, $a_2=\min v_2$, and $a_3=\sigma(\tau)$. By definition, we know that $a_2<\min V(\tau)$, so that in the diagram representation of the three vertices $v_1,v_2,V(\tau)$ (Figure~\ref{fig:diaggluing}(a)) we know that the blocks $A_1,A_2,A_3,B_2,B_3$ are empty. Moreover, since $\tau$ is a trisection, $a_3$ is not the minimum inside its vertex, so the block $C_3$ is \emph{not} empty. Hence, comparing Figures~\ref{fig:diaggluing}(c) and~\ref{fig:trisectionmarquee}, and observing once again that the blocks $C_3$ and $K$ coincide, we see that after the gluing, $\tau$ is a trisection of type {\bf II} in the new map $\barM$. Moreover, since the slicing and gluing operations are inverse one to the other, it is clear that $\Gamma(\barM,\tau)=(\M,v_1,v_2,\tau)$. Hence, $\Gamma\circ\Psi$ is the identity, and the proposition is proved.
\end{proof}

\section{Iterating the bijection.}

Of course Proposition~\ref{prop:type1} looks nicer than its counterpart Proposition~\ref{prop:type2}: in the first one, one only asks to distinguish three vertices in a map of lower genus, whereas in the second one, the marked triple must satisfy a nontrivial constraint (Equation~\ref{eq:contrainteordre}). In this section we will work a little more in order to get rid of this problem. We start with two definitions (observe that for $k=3$ this is coherent with what precedes):
\begin{definition}
We let $\mathcal{U}^k_g(n)$ be the set of unicellular maps of genus $g$ with $n$ edges, and $k$ distinct unordered distinguished vertices.
\end{definition}
\begin{definition}
We let $\mathcal{D}_g(n)=\mathcal{D}^I_g(n)\cup\mathcal{D}^{II}_g(n)$ be the set of unicellular maps of genus $g$ with $n$ edges, and a distinguished trisection.
\end{definition}

\subsection{Training examples: genera $1$ and $2$.}

Observe that the set $\mathcal{V}_0(n)$ is empty, since there are no trisections in a plane tree. Hence, from Proposition~\ref{prop:type2}, there are no trisections of type {\bf II} in a map of genus~$1$ (i.e. $\mathcal{D}_1^{II}(n)=\varnothing$). Proposition~\ref{prop:type1} gives:
\begin{corollary}
The set $D_1(n)$ of unicellular maps of genus $1$ with $n$ edges and a distinguished trisection is in bijection with the set\/ $U_0^3(n)$ of rooted plane trees with $n$ edges and three distinguished vertices. 
\end{corollary}
\noindent Since from the trisection lemma (Lemma~\ref{lemma:trisection}) each unicellular map of genus $1$ has exactly $2$ trisections, we obtain that the number $\epsilon_1(n)$ of rooted unicellular maps of genus $1$ with $n$ edges satisfies:
$$
2 \cdot \epsilon_1(n) = {n+1 \choose 3} \mathrm{Cat}(n),
$$
which gives a neat combinatorial proof of the formula $\epsilon_1(n)=\frac{(n+1)n(n-1)}{12}\mathrm{Cat}(n)$~\cite{Lehman-Walsh-genus-I}.\\

We now consider the case of genus $2$. Let $\M$ be a unicellular map of genus $2$, and $\tau$ be a trisection of $\M$. If $\tau$ is of type {\bf I}, we know that we can use the application $\Phi^{-1}$, and obtain a unicellular map of genus~$1$, with three distinguised vertices.

Similarly, if $\tau$ is of type {\bf II}, we can apply the mapping $\Psi^{-1}$, and we are left with a unicellular map $\M'$ of genus $1$, and a marked triple $(v_1,v_2,\tau)$, such that $\min_{\scriptscriptstyle \M'} v_1 <_{\scriptscriptstyle \M'} \min_{\scriptscriptstyle \M'} v_2 <_{\scriptscriptstyle \M'} \min_{\scriptscriptstyle \M'} V(\tau)$. From now on, we use the more compact notation: $v_1 <_{\scriptscriptstyle \M'} v_2 <_{\scriptscriptstyle \M'} V(\tau)$, i.e. we do not write the min's anymore.
The map $(\M',\tau)$ is a unicellular map of genus $1$ with a distinguished trisection: therefore we can apply the mapping $\Phi^{-1}$ to $(\M',\tau)$. We obtain a plane tree $\M''$, with three distinguished vertices $v_3,v_4,v_5$ inherited from the slicing of $\tau$ in $\M'$; since  those three vertices are undistinguishable, we can assume that $ v_3<_{\scriptscriptstyle \M''} v_4 <_{\scriptscriptstyle \M''} v_5$.
Observe that in $\M''$ we also have the two marked vertices $v_1$ and $v_2$ inherited from the slicing of $\tau$ in $\M$. Moreover the fact that $ v_1 <_{\scriptscriptstyle \M'}  v_2 <_{\scriptscriptstyle \M'} V(\tau)$ in $\M'$ implies that $v_1 <_{\scriptscriptstyle \M''} v_2 <_{\scriptscriptstyle \M''} v_3$ in $\M''$: this follows from Remark~\ref{rem:ordre}. 
%
%
Hence, we are left with a plane tree $\M''$, with \emph{five} distinguished vertices 
$v_1<_{\scriptscriptstyle \M''}v_2<_{\scriptscriptstyle \M''}v_3<_{\scriptscriptstyle \M''}v_4<_{\scriptscriptstyle \M''}v_5$. 
Conversely, given such a $5$-tuple of vertices, it is always possible to glue the three last ones together by the mapping $\Phi$ to obtain a triple $(v_1,v_2,\tau)$ satisfying Equation~\ref{eq:contrainteordre}, and then to apply the mapping $\Psi$ to retrieve a map of genus $2$ with a marked trisection of type {\bf II}. This gives:
\begin{corollary}
The set $\mathcal{D}^{II}_2(n)$ is in bijection with the set $\mathcal{U}^5_0(n)$ of plane trees with \emph{five} distinguished vertices.

The set $\mathcal{D}_2(n)$ of unicellular maps of genus $2$ with one marked trisection is in bijection with the set $\mathcal{U}_1^3(n)\cup\mathcal{U}_0^5(n)$.
\end{corollary}
\noindent From Euler's formula, a unicellular map of genus $1$ with $n$ edges has $n-1$ vertices.
Since from the trisection lemma each unicellular maps of genus $2$ has $4$ trisections,  we obtain the following formula for the number $\epsilon_2(n)$ of unicellular maps of genus $2$ with $n$ edges:
$$
4\cdot \epsilon_2(n) = {n-1 \choose 3} \epsilon_1(n) + {n+1 \choose 5} \mathrm{Cat}(n),
$$
from which it follows that 
$$
\epsilon_2(n)= \frac{(n+1)n(n-1)(n-2)(n-3)(5n-2)}{1440} \mathrm{Cat}(n).
$$

\subsection{The general case, and our main theorem.}

In the general case, we will work as in the example of genus $2$: while we find trisections of type {\bf II}, we open them (and we decrease the genus of the map), and we stop when we have opened the first encountered trisection of type {\bf I}. We start with the description of the inverse procedure, which goes as follows.\\

We let $p\geq 0$ and $q\geq 1$ be two integers, and $(\M,v_*)=(\M,v_1,\dots,v_{2q+1})$ be an element of $\mathcal{U}^{2q+1}_p(n)$. Up to renumbering the vertices, we can assume that $v_1 <_{\scriptscriptstyle \M} v_2 <_{\scriptscriptstyle \M} \dots <_{\scriptscriptstyle \M} v_{2q+1}$.
\begin{definition}
We consider the following procedure:
\\ \noindent {\bf i.} 
Glue the three last vertices $v_{2q-1},v_{2q},v_{2q+1}$ together, via the mapping $\Phi$, in order to obtain a new map $\M_{1}$ of genus $p+1$ with a distinguished trisection $\tau$ of type {\bf I}.
\\ \noindent {\bf ii.} 
{\bf for} $i$ from $1$ to $q-1$ {\bf do:}\\
Let $(v_{2q-2i-1},v_{2q-2i},\tau)$ be the triple consisting of the last two vertices which have not been used until now, and the trisection $\tau$. Apply the mapping $\Psi$ to that triple, in order to obtain a new map $\M_{i+1}$ of genus $p+i+1$, with a distinguished trisection $\tau$ of type {\bf II}.
\\ \noindent {\bf end for.}\\
 We let $\Lambda(\M,v_*):=(\M_{q},\tau)$ be the map with a distinguished trisection obtained at the end of this procedure. Observe that if $q=1$, the distinguished trisection is of type {\bf I}, and that it is of type {\bf II} otherwise.
\end{definition}

As in the case of genus $2$, we have the next theorem:
\begin{theorem}[Our main result]\label{thm:main}
The application $\Lambda$ defines a bijection:
$$
\Lambda : \  \biguplus _{p=0}^{g-1} \mathcal{U}^{2g-2p+1}_p(n)
\longrightarrow \mathcal{D}_g(n).
$$
In other words, all unicellular maps of genus $g$ with a distinguished trisection can be obtained in a canonical way by starting with a map of a lower genus with an odd number of distinguished vertices, and then applying once the mapping $\Phi$, and a certain number of times the mapping $\Psi$.
\end{theorem}

Given a map with a marked trisection $(\M,\tau)$, the converse application consists in slicing recursively the trisection $\tau$ while it is of type {\bf II}, then slicing \emph{once} the obtained trisection of type {\bf I}, and remembering all the vertices resulting from the successive slicings.
More formally, we have the next proposition:
\begin{proposition}\label{prop:main}
Let $(\M,\tau)\in\mathcal{D}_g(n)$. We define the pair $\Xi(\M,\tau)$ by the following procedure:
\begin{enumerate}
\setlength{\itemsep}{-2pt}
\item We let $\M_0:=\M$  and $i:=0$.
\item \label{goto} 
{\bf If} $\tau$ is of type {\bf II} in $\M_i$, we let
$(\M_{i+1},v_{2i+1},v_{2i+2}):=\Psi^{-1}(\M_i,\tau)$. Then we let $i:=i+1$ and we return to step~\ref{goto}. \\
{\bf  Else,} $\tau$ is of type {\bf I} in $\M_i$ and we go to step~\ref{goto1}.
\item \label{goto1} Let $(\M_{i+1},v_{2i+1},v_{2i+2},v_{2i+3}):=\Phi^{-1}(\M_i,\tau)$.
We end the pocedure and we let 
$$\Xi(\M,\tau):=(\M_{i+1}, v_1,v_2,\dots,v_{2i+3}).$$
\end{enumerate}
\noindent
Then the mapping
$$
\Xi : \  \mathcal{D}_g(n)
\longrightarrow \biguplus _{p=0}^{g-1} \mathcal{U}^{2g-2p+1}_p(n).
$$
is a bijection, which is the inverse bijection of $\Lambda$.
\end{proposition}

\begin{proof}[Proof of Theorem~\ref{thm:main} and Proposition~\ref{prop:main}]
First, the mapping $\Xi$ is well defined. Indeed, by definition of a trisection of type {\bf II}, we know by induction that each time we enter steps~\ref{goto} and~\ref{goto1}, $\tau$ is a trisection of the map $\M_i$. Moreover, since the genus of the maps $\M_i$ decreases with $i$, we necessarily reach step~\ref{goto1}, and the procedure stops.

Then, the mapping $\Lambda$ is clearly injective, since the mappings $\Psi$ and $\Phi$ are.

Finally, to prove at the same time that $\Xi$ is injective and that it is the inverse mapping of $\Lambda$, it is enough to show that the vertices $v_i$ produced by the procedure defining $\Xi$ satisfy $v_1 \infM v_2 \infM \dots \infM v_{2q+1}$. Indeed, after that it will be clear by construction that $\Lambda\circ\Xi=\Xi\circ\Lambda=Id$.
Now, we deduce from Remark~\ref{rem:ordre} and by an induction on $i$ that after the $i$th passage in step~\ref{goto} in the definition of $\Xi$, we have
$v_1 \infMii v_2 \infMii \dots \infMii v_{2i+2}$. The same remark shows that the end of step~\ref{goto1}, we have
$v_1 \infMii v_2 \infMii \dots \infMii v_{2i+3}$, which concludes the proof.
\end{proof}

\section{Enumerative corollaries.}

\subsection{A combinatorial identity}

Using the trisection lemma (Lemma~\ref{lemma:trisection}), Euler's formula, and Theorem~\ref{thm:main}, we obtain the following new identity (stated in the introduction as Equation~\ref{eq:identity}): 
\begin{theorem}\label{thm:identity}
The number $\epsilon_g(n)$ of rooted unicellular maps of genus $g$ with $n$ edges satisfies the following combinatorial identity:
\begin{eqnarray*}
 2g  \cdot \epsilon_g(n)\!\!\!\!\! &=& \!\!\!\!\! {n+3-2g\choose 3} \epsilon_{g-1}(n) + {n+5-2g\choose 5} \epsilon_{g-2}(n) + \dots + {n+1\choose 2g+1} \epsilon_{0}(n).
\end{eqnarray*}
\end{theorem}
\noindent Observe that this identity is recursive only in the genus (the number of edges $n$ is fixed). For that reason, it enables one to compute easily, for a fixed $g$, the closed formula giving $\epsilon_g(n)$ by a simple iteration (as we did for genera $1$ and $2$).

\subsection{The polynomials $R_g(n)$}

By iterating Theorem~\ref{thm:main}, we obtain that all unicellular maps of genus $g$ with $n$ edges can be obtained from a plane tree with $n$ edges, by successively gluing vertices together. From the enumeration viewpoint, we obtain the first bijective proof that the numbers $\epsilon_g(n)$ are the product of a polynomial and a Catalan number:
\begin{corollary}
The number $\epsilon_g(n)$ of unicellular maps of genus $g$ with $n$ edges equals:
$$
\epsilon_g(n)=R_g(n)\mathrm{Cat}(n),
$$
where $R_g$ is the polynomial of degree $3g$ defined by the formula:
$$
R_g(n) = \sum_{0=g_0<g_1<\dots<g_r=g}\,\,\prod_{i=1}^r \frac{1}{2g_i}{{n+1-2g_{i-1}}\choose{2(g_i-g_{i-1})+1}}.
$$
\end{corollary}
\begin{proof}
The statement directly comes from an iteration of the bijection of Theorem~\ref{thm:main}. More precisely, the formula given here for $R_g(n)$ reads as follows. In order to generate a unicellular map of genus~$g$, we start with a plane tree with $n$ edges, and we apply a certain number of times (say $r$) the mapping~$\Lambda$ to create unicellular maps of increasing genera. In the formula, $g_1< \dots < g_r = g$ are the genera of the maps produced by the successive applications of~$\Lambda$. Now, in order to increase the genus from $g_{i-1}$ to $g_i$, we have to choose $2(g_i-g_{i-1})+1$ vertices in a unicellular map of genus $g_{i-1}$, which gives the binomial in the formula. The factor $1/(2g_i)$ is here to compensate the multiplicity in the construction coming from the trisection lemma (Lemma~\ref{lemma:trisection}).
\end{proof}
\noindent From Theorem~\ref{thm:identity} and the fact that $\mathrm{Cat}(n)$ is asymptotically equivalent to $\frac{1}{\sqrt{\pi}}n^{-\frac{3}{2}}4^n$, one recovers easily the asymptotic behaviour of $\epsilon_g(n)$:
\begin{corollary}
The polynomial $R_g(n)$ has degree $3g$ and leading coefficient $\tfrac{1}{12^g g!}$. When $n$ tends to infinity, one 
has~\cite{BeCaRo:asymptotic-nb-tree-rooted}:
$$
\epsilon_g(n)\sim \frac{1}{12^gg! \sqrt{\pi}} n^{3g-\frac{3}{2}} 4^n.
$$
\end{corollary}

\noindent Our construction also answers a question of Zagier concerning the interpretation of a property of the polynomials $R_g$:
\begin{corollary}[Zagier {\cite[p. 160]{LaZv}}]
For each integer $g\geq 1$, the polynomial $R_g(n)$ is divisible by $(n+1) \dots (n+1-2g)$.
\end{corollary}
\begin{proof}
For any fixed sequence $g_1< \dots < g_r$ of intermediate genera, one has to choose $\sum_{i}(2g_i-2g_{i-1}+1)=2g+r$ distinct vertices in the original plane tree in order to apply the construction. Since $r\geq 1$, this number is at least $2g+1$, and the statement follows.
\end{proof}

\section{Variants.}

\subsection{Bipartite unicellular maps}

A unicellular map is \emph{bipartite} if one can color its vertices in black and white in such a way that only vertices of different colors are linked by an edge. By convention, the root-vertex will always be colored in white.
\begin{definition}
We let $\beta_g(i,j)$ be the number of bipartite unicellular maps of genus $g$ with $i$ white vertices and $j$ black vertices. Such a map has $i+j+2g-1$ edges.
\end{definition}
It is clear that that our construction also applies to bipartite unicellular maps: the only difference is that, if we want the gluing operations $\Phi$ and $\Psi$ to preserve the bipartition of the map, we have to paste together only vertices \emph{of the same color}. We therefore obtain the following variant of our main identity:
\begin{proposition}The number $\beta_g(i,j)$ of bipartite unicellular maps with $i$ white vertices and $j$ black vertices obey the following recursion formula:
\begin{eqnarray}
\label{eq:identitebip}
2g \cdot \beta_g(i,j) = 
\sum_{p=0}^{g-1} \binom{i+2g-2p}{2g-2p+1} \beta_{p}(i+2g-2p,j)
\ +\ 
\sum_{p=0}^{g-1} \binom{j+2g-2p}{2g-2p+1} \beta_{p}(i,j+2g-2p).
\end{eqnarray}
\end{proposition}

\begin{corollary}
We have $\beta_g(i,j)=S_g(i,j) \beta_0(i,j)$, where $\beta_0(i,j)=\frac{i+j-1}{ij}\binom{i+j-2}{i-1}^2$ is the number of bipartite plane trees with $i$ white vertices and $j$ black vertices computed in~\cite{Goulden-Jackson:book}, and $S_g$ the polynomial in $(i,j)$ defined by:
$$
S_g(i,j) = \frac{1}{2g} \sum_{p=0}^{g-1} \binom{i+2g-2p}{2g-2p+1} S_{p}(i+2g-2p,j)
\ +\ 
 \frac{1}{2g}\sum_{p=0}^{g-1} \binom{j+2g-2p}{2g-2p+1} S_{p}(i,j+2g-2p)
$$
with $S_0=1$.
\end{corollary}
\noindent For example for the first genera we obtain:
\begin{align*}
S_0=1\ \ ; \ \ S_1(i,j)=\frac{(i+2)(i+1)i+(j+2)(j+1)j}{12} \ \ ; \ \ S_2(i,j) = s_2(i,j)+ s_2(j,i)
\end{align*}
where
$s_2(i,j)=$ {\small $\frac{i (i+1) (i+2) (i^5+22 i^4+211 i^3+2 i^2 j+998 i^2+i^2 j^3+3 i^2 j^2+21 i j^2+2248 i+7 i j^3+14 i j+96 j^2+1920+64 j+32 j^3)
}{5760} $}.\\

\subsection{Precubic unicellular maps}

A unicellular map is \emph{precubic} if all its vertices have degree $1$ or $3$. In such a map, all trisections are necessarily of type {\bf I}: indeed, a trisection of type {\bf II} cannot appear in a vertex of degree less than $4$. Therefore, each precubic map can be obtained in exactly $2g$ different ways from a \emph{precubic} map of genus $(g-1)$ with three distinguished \emph{leaves}. By repeating the argument $g$ times, we see that each precubic unicellular map can be otbained in exactly $2g \cdot 2(g-1) \dots 2 = 2^gg!$ different ways from a precubic tree (a plane tree where all vertices have degrees $1$ or $3$), by $g$ successive gluings of triple of leaves.

Now, we can easily enumerate precubic trees with $n$ edges. First, we observe that by removing a leaf to such a tree, we find a binary tree with $n-1$ edges (and $n$ vertices). This implies that $n=2m+1$, where $m$ is the number of \emph{nodes} of the binary tree, and that the number of precubic trees with $n$ edges which are rooted on a leaf is the Catalan number $\mathrm{Cat}(m)$. A double-counting argument then shows that those whose root-vertex has degree $3$ are counted by the number $\frac{3m}{m+2}\mathrm{Cat}(m)$: indeed, the number $3m\mathrm{Cat}(m)$ counts precubic trees which are rooted at the same time on a leaf and a vertex of degree $3$, and these trees can also be obtained by distinguishing one of the $(m+2)$ leaves in a tree which is rooted on a vertex of degree $3$. Finally, the number of all precubic rooted trees with $n$ edges equals $\displaystyle\big(1+\frac{3m}{m+2}\big)\Cat(m)$. We thereby obtain:
\begin{corollary}
The number $\xi_g(n)$ of precubic unicellular maps of genus $g$ with $n=2m+1$ edges is given by:
$$
\xi_g(n) = \frac{1}{2^gg!}{m+2 \choose \ 3,3,\dots , 3, m+2-3g\ }\Big(1+\frac{3m}{m+2} \Big) \Cat(m) = \frac{(4m+2)(2m)!}{12^gg!(m+2-3g)!m!}.
$$
\end{corollary}

\noindent Precubic unicellular maps which have no leaves necessarily have $6g-3$ edges. These objects appeared previously in the litterature~(\cite{Lehman-Walsh-genus-I, Bacher-Vdovina}, and recently in~\cite{ChMaSc} under the name of \emph{dominant schemes}). We can recover their number by setting $m=3g-2$ in the previous equation. In that case, the bijection given here reduces to the one given in our older paper~\cite{Chapuy:trisections}, in which the following corollary already appeared. However, we copy it out here for completeness:
\begin{corollary}[\cite{Lehman-Walsh-genus-I}]\label{cor:precubiques}
The number of rooted unicellular maps of genus $g$ with all vertices of degree~$3$ is:
$$
\frac{2(6g-3)!}{12^g g!(3g-2)!}.
$$
Dually, this number counts rooted triangulations of genus $g$ with only one vertex.
\end{corollary}

\subsection{Labelled unicellular maps.}

A \emph{labelled unicellular map} is a pair $(\M,l)$ such that $\M$ is a rooted unicellular map, and $l$ is a \emph{labelling} of the vertices of $\M$, i.e. a mapping 
$l: V(\M)\rightarrow \mathbb{Z}$ such that if $v_1, v_2$ are two adjacent vertices in~$\M$, then $l(v_1)-l(v_2)\in\{-1,0,1\}$, and such that the root-vertex has label $0$. Our interest for these objects comes from the so-called Marcus-Schaeffer bijection:
\begin{proposition}[Marcus Schaeffer \cite{MaSc}, see also \cite{ChMaSc} for the version stated here]
Let $m_g(n)$ be the number of (all) rooted maps of genus $g$ with $n$ edges, and let $l_g(n)$ be the number of labelled unicellular maps of genus $g$ with $n$ edges. Then one has:
$$
(n+2)\, m_g(n) = 2 \, l_g(n).
$$
This identity comes from an explicit and simple bijection.
\end{proposition}
Therefore it is interesting to see what our construction becomes on these labelled unicellular maps. We let $\mathcal{L}^{(k)}_g(n)$ be the set of rooted labelled unicellular maps of genus $g$ with $n$ edges and $k$ distinguished vertices \emph{carrying the same label}. We also let $\mathcal{DL}_g(n)$ be the set of labelled unicellular maps carrying a distinguished trisection. We have:
\begin{corollary}
The application $\Lambda$ induces a bijection:
$$
\Lambda : \  \biguplus_{p=0}^{g-1} \mathcal{L}^{(2g-2p+1)}_p(n)
\longrightarrow \mathcal{DL}_g(n).
$$
\end{corollary}
\begin{proof}
The only thing to change in our construction so that it is well-defined is to restrict the gluing operation to vertices of the same label, which is exactly what we do here. 
\end{proof}
Observe that it is not obvious to compute the cardinality of $\mathcal{L}^{(k)}_g(n)$: in order to compute it from $l_g(n)$, one would need non trivial information about the repartition of labels of vertices in a randomly labelled unicellular map of genus $g$, or by induction, in a randomly labelled \emph{plane tree} (see~\cite{Chapuy:trisections} for the connection with functionals of continuum random trees, in the case of $k=3$).
Therefore, in order to transform the last corollary into a counting strategy for maps on surfaces, one would need to understand much better the structure of plane labelled trees than we do at the moment. 

\section{More computations.}\label{subsec:guitter}

We now sketch a computation of Emmanuel Guitter~\cite{Guitter:comm}, that enables to recover the Harer-Zagier formula from our construction.
For all $n\geq 1$, we let 
$
F_n(x) = \sum_{g\geq 0}\epsilon_g(n) x^{n+1-2g}
$
be the generating function of unicellular maps with $n$ edges, where the variable $x$ marks the number of vertices. Then we have:
\begin{proposition}[\cite{Guitter:comm}]\label{prop:Guittereq}
For every real sequence $(a_n)_{ n\geq 0}$, the formal power series $F(x,y)=\sum_{n\geq 0} a_n y^{n+1} F_n(x)$ is solution of the differential equation:
\begin{eqnarray}\label{eq:Guittereq}
2y \cdot\frac{\partial }{\partial y} F(x,y) = F(x+1,y)-F(x-1,y)
\end{eqnarray}
\end{proposition}
\begin{proof}
Clearly, the series $\displaystyle \frac{1}{2}\big(F_n(x+1)-F_n(x-1)\big)$ is the generating function of unicellular maps with $n$ edges, in which an \emph{odd} number of vertices have been distinguished, and are no longer counted in the exponant of $x$. These objects are divided into two categories: either the number of distinguished vertices is $\geq 3$, or it is equal to one. By our main theorem, objects of the first kind are in bijection with unicellular maps with $n$ edges and a distinguished trisection; objects of the second kind are unicellular maps with $n$ edges with a distinguished vertex. Now, by the trisection lemma and Euler's formula, in each map the number of trisections plus the number of vertices equals $n+1$. Therefore we have:
$\displaystyle \frac{1}{2}\big(F_n(x+1)-F_n(x-1)\big) = (n+1) F_n(x)$ and the proposition follows.
\end{proof}
\begin{corollary}[\cite{Harer-Zagier}]
Let $a_n=\frac{2^n n!}{(2n)!}$, and let $F(x,y)=\sum_{n\geq 0} a_n y^{n+1} F_n(x)$. Then one has:
\begin{eqnarray}\label{eq:IZ}
F(x,y) = \frac{1}{2}\left(\frac{1+y}{1-y}\right)^x -\frac{1}{2}.
\end{eqnarray}
\end{corollary}
Note that the Harer-Zagier formula (Equation~\ref{eq:HZformula}) can be recovered by expanding the $x$-th power in~(\ref{eq:IZ}).
\begin{proof} We follow~\cite{Guitter:comm}. First, 
one easily checks  that the function given here is solution of Equation~\ref{eq:Guittereq}. Moreover, a solution to Equation~\ref{eq:Guittereq} is entirely characterized by its "planar terms", i.e. by the coefficients of $x^{n+1}y^{n+1}$ for all $n\geq 0$ (think about computing the coefficients inductively via Equation~\ref{eq:identity}). Hence the only thing to check is that $\lim_{y\rightarrow 0 }F(\frac{u}{y},y)$ is equal to $\sum_{n\geq 0} a_n \mathrm{Cat}(n)u^{n+1} = \frac{1}{2}\big(\exp(2u)-1\big)$, which is immediate from~(\ref{eq:IZ}).
\end{proof}

\noindent We conclude with an extension of the previous computation to bipartite unicellular maps. For these maps, the ordinary generating series do not have a closed formula~\cite{Adrianov}, and it is better to work with "modified" generating series. More precisely, following Adrianov~\cite{Adrianov}, we introduce for each integer $v \geq 0$ the series $\phi_v(x):=\sum_{k\geq 1} k^v x^{k-1}$. We consider the \emph{modified generating series} of bipartite unicellular maps defined as follows:
$$ 
B_n(x,y) = \sum_{i,j\geq 0} \beta_{\frac{n+1-i-j}{2}}(i,j) \phi_i(x)\phi_j(y).
$$
Observe that by Euler's formula, a unicellular map with $i+j$ vertices and genus $\frac{n+1-i-j}{2}$ has $n$ edges, so that
$B_n(x,y)$ is the generating function of bipartite unicellular maps with $n$ edges, in which a map with $i$ white and $j$ black vertices is given a weight $\phi_i(x)\phi_j(y)$.  By investigating the effect of the deletion of vertices in the context of modified generating series, one obtains the following analogue of Proposition~\ref{prop:Guittereq}.
\begin{proposition}
The formal power series $B_n(x,y)$ is solution of the differential equation:
\begin{eqnarray}\label{eq:diffbiparties}
  \frac{\partial }{\partial x}\big( (1-x^2)B_n (x,y) \big)
+ \frac{\partial }{\partial y}\big( (1-y^2)B_n (x,y) \big)
= (n+1) B_n(x,y).
\end{eqnarray}
\end{proposition}
\begin{corollary}[\cite{Adrianov}]
The series $B_n(x,y)$ admits the following closed form:
\begin{eqnarray}\label{eq:adrianov}
B_n(x,y)=n! \frac{(1-xy)^{n-1}}{((1-x)(1-y))^{n+1}}.
\end{eqnarray}
\end{corollary}
\begin{proof}
One easily checks that the series given here is solution of Equation~\ref{eq:diffbiparties}.
Now, as in the monochromatic case, a formal power series 
$$C_n(x,y)= \sum_{1\leq i+j\leq n+1} c_{\frac{n+1-i-j}{2}}(i,j) \phi_i(x)\phi_j(y)$$ 
which is a solution of Equation~\ref{eq:diffbiparties} is characterized by its "planar terms", i.e. by the sequence of numbers $(c_{0}(i,j))_{i\geq1, j\geq 1}$. 
Therefore it is enough to prove that the numbers $c_0(i,j)$ corresponding to the function $C_n(x,y):=n! \frac{(1-xy)^{n-1}}{((1-x)(1-y))^{n+1}}$ are equal to the numbers $\beta_0(i,j)=\frac{i+j-1}{ij}\binom{i+j-2}{i-1}^2$.

Now, set $X=\frac{1}{1-x}$, $Y=\frac{1}{1-y}$ so that $C_n(x,y)=n!X^2Y^2(X+Y-1)^{n-1}$ is a polynomial $\widetilde{C}_n(X,Y)$ in $X$ and $Y$.
Using the fact that around $x=1$, one has $\phi_v(x)=\tfrac{v!}{\ \ (1-x)^{v+1}} + O\left(\tfrac{1}{(1-x)^{v}}\right)$, one obtains that
for all $i,j$ such that $i+j=n+1$,  the coefficient of $X^{i+1}Y^{j+1}$ in the polynomial $\widetilde{C}_n(X,Y)$ is $i!j!c_{0}(i,j)$.
Therefore we have:
$$c_0(i,j)= \frac{n!}{i!j!} \times \mathrm{Coeff}_{X^{i-1}Y^{j-1}} \Big( (X+Y-1)^{n-1} \Big)= \frac{n!}{i!j!} {n-1 \choose i-1,j-1,0} = \beta_0(i,j).$$
\end{proof}

\bibliographystyle{alpha}
\bibliography{UnicellularExact}
\end{document}